\documentclass{amsart}

\usepackage{amsmath}
\usepackage{amsfonts}
\usepackage{amssymb}
\usepackage{color}
\usepackage{pstricks,pst-plot,pstricks-add,pst-math}
\usepackage{enumerate}
\usepackage{verbatim}
\usepackage{epsfig}
\usepackage{setspace}

\IfFileExists{myowntimes.sty}{\usepackage{myowntimes}\usepackage{mathrsfs}}
        {
\usepackage{mathrsfs}}

\DeclareFontFamily{OT1}{eusb}{} \DeclareFontShape{OT1}{eusb}{m}{n}
{<5> <6> <7> <8> <9> <10> <11> <12> <14.4> eusb10}{}
\DeclareMathAlphabet{\eusb}{OT1}{eusb}{m}{n}

\DeclareFontFamily{OT1}{eusm}{} \DeclareFontShape{OT1}{eusm}{m}{n}
{<5> <6> <7> <8> <9> <10> <11> <12> <14.4> eusm10}{}
\DeclareMathAlphabet{\eusm}{OT1}{eusm}{m}{n}

\DeclareFontFamily{OT1}{eufm}{} \DeclareFontShape{OT1}{eufm}{m}{n}
{<5> <6> <7> <8> <9> <10> <11> <12> <14.4> eufm10}{}
\DeclareMathAlphabet{\mathfrak}{OT1}{eufm}{m}{n}

\DeclareFontFamily{OT1}{fraktura}{}
\DeclareFontShape{OT1}{fraktura}{m}{n} {<5> <6> <7> <8> <9> <10> <11>
  <12> <13> <14.4> [1.1] eufm10}{}
\DeclareMathAlphabet{\fraktura}{OT1}{fraktura}{m}{n}

\DeclareFontFamily{OT1}{cmfi}{} \DeclareFontShape{OT1}{cmfi}{m}{n}
{<5> <6> <7> <8> <9> <10> <11> <12> <13> <14.4> [0.9] cmfi10}{}
\DeclareMathAlphabet{\cmfi}{OT1}{cmfi}{b}{n}

\DeclareFontFamily{OT1}{cmss}{} \DeclareFontShape{OT1}{cmss}{m}{n}
{<5> <6> <7> <8> <9> <10> <11> <12> <13> <14.4> cmss10}{}
\DeclareMathAlphabet{\cmss}{OT1}{cmss}{m}{n}

\theoremstyle{plain} \newtheorem{Lemma}{Lemma}[section]
\theoremstyle{plain} \newtheorem{lemma}{Lemma}[section]
\theoremstyle{plain} 
\theoremstyle{plain} \newtheorem{proposition}[Lemma]{Proposition}

\theoremstyle{plain} 
\theoremstyle{plain} 
\theoremstyle{plain} \newtheorem{theorem}[Lemma]{Theorem}

\theoremstyle{plain} 
\theoremstyle{plain}

\theoremstyle{definition} 
\theoremstyle{definition} \newtheorem{definition}[Lemma]{Definition}

\theoremstyle{remark} \newtheorem{remark}[Lemma]{Remark}
\theoremstyle{remark}

\numberwithin{equation}{section}

\renewcommand{\section}{\secdef\sct\sect}
\newcommand{\sct}[2][default]{%
\refstepcounter{section}
\addcontentsline{toc}{section}{{\tocsection
    {}{\thesection}{\!\!\!\!#1\dotfill}}{}}
\vspace{0.7cm}
\centerline{\scshape\thesection.\ #1} \nopagebreak \vspace{0.2cm}}
\newcommand{\sect}[1]{%
\vspace{0.4cm} \centerline{\large\scshape\rmfamily #1}
\vspace{0.2cm}
}

\renewcommand{\subsection}{\secdef\subsct\sbsect}
\newcommand{\subsct}[2][default]{\refstepcounter{subsection}
\addcontentsline{toc}{subsection}
{{\tocsection{\!\!}{\hspace{1.2em}\thesubsection}{\!\!\!\!#1\dotfill}}{}}
\nopagebreak
{\flushleft\bf
\thesubsection~\bf #1.~}
\noindent
\nopagebreak}
\newcommand{\sbsect}[1]{
\noindent
\textbf{#1.~}
}

\renewcommand{\subsubsection}{%
\secdef \subsubsect\sbsbsect}
\newcommand{\subsubsect}[2][default]{%
\refstepcounter{subsubsection}
\addcontentsline{toc}{subsubsection}{{\tocsection{\!\!}
{\hspace{3.05em}\thesubsubsection}{\!\!\!\!#1\dotfill}}{}}
\nopagebreak
\vspace{0.15\baselineskip} \nopagebreak {\flushleft\rmfamily
\itshape\thesubsubsection
\ \rmfamily #1\/.}\ }
\newcommand{\sbsbsect}[1]{\vspace{0.1cm}\noindent
\rmfamily \itshape
\arabic{section}.\arabic{subsection}.\arabic{subsubsection} \
\sffamily #1\/.\ }

\renewcommand{\caption}[1]{%
\vglue0.5cm
\refstepcounter{figure}
\begin{minipage}{0.9\textwidth}\small {\sc Figure~\thefigure. }#1\end{minipage}}

 \newcommand{\ie}{\emph{i.e.,}}

 \newcommand{\iid}{i.i.d.}

 \DeclareMathOperator{\E}{\mathbb{E}}
 \DeclareMathOperator{\pr}{\mathbb{P}}

 \newcommand{\mC}{\mathcal {C}}

 \newcommand{\mP}{\mathcal{P}}

 \newcommand{\mV}{\mathcal{V}}
 \newcommand{\mW}{\mathcal{W}}

 \newcommand{\bP}{\mathbb{P}}
 
 \newcommand{\bR}{\mathbb{R}}

 \newcommand{\ga}{\alpha}

 \newcommand{\gC}{\Gamma}
 \newcommand{\gd}{\delta}

 \newcommand{\gs}{\sigma}

\newtheorem{thm}{Theorem}

\newtheorem{corollary}[thm]{Corollary}

\newcommand{\BB}{\mathcal B}
\newcommand{\CC}{\mathcal C}

\newcommand{\PP}{\mathcal P}

\newcommand{\BbbP}{\mathbb P}

\def\myffrac#1#2 in #3{\raise 2.6pt\hbox{$#3 #1$}\mkern-1.5mu\raise 0.8pt\hbox{$#3/$}\mkern-1.1mu\lower 1.5pt\hbox{$#3 #2$}}
\newcommand{\ffrac}[2]{\mathchoice%
        {\myffrac{#1}{#2} in \scriptstyle}
        {\myffrac{#1}{#2} in \scriptstyle}
        {\myffrac{#1}{#2} in \scriptscriptstyle}
        {\myffrac{#1}{#2} in \scriptscriptstyle}
}

\begin{document}

\title[Scaling Limits for Posets]{Scaling Limits for Width Two
  Partially Ordered Sets: The Incomparability Window}

\author{Nayantara Bhatnagar$^*$}
\thanks{$^*$ Department of Statistics, UC Berkeley, {\tt
    email:nayan@stat.berkeley.edu}, Supported by DOD ONR grant
  N0014-07-1-05-06 and DMS 0528488.}

\author{Nick Crawford$^\dagger$}
\thanks{$^\dagger$ Department of Statistics, UC Berkeley, {\tt
    email:crawford@stat.berkeley.edu}, Supported by DMS 0548249
  (CAREER)}

\author{Elchanan Mossel$^\ddagger$}
\thanks{$^\ddagger$ Department of Statistics and Department of Computer
  Science , UC Berkeley, {\tt email:mossel@stat.berkeley.edu},
  Supported by DOD ONR grant
  N0014-07-1-05-06, DMS 0528488, DMS 0548249 (CAREER), and a Sloan
  Fellowship.}

\author{Arnab Sen$^\S$}
\thanks{$^\S$ Department of Statistics, UC Berkeley, {\tt
    email:arnab@stat.berkeley.edu}, Supported by DOD ONR grant
  N0014-07-1-05-06, DMS 0528488, and DMS 0548249 (CAREER)}

\begin{abstract}
We study the structure of a uniformly randomly chosen partial
order of width 2 on $n$ elements. 
We show that under the appropriate
scaling, the number of incomparable elements converges to the height
of a one dimensional Brownian excursion at a uniformly chosen random
time in the interval $[0,1]$, which follows the Rayleigh
distribution.
\end{abstract}

\maketitle

\section{Introduction}
The study of the typical local structure of large random combinatorial
objects is of central interest  
in modern combinatorics. Classical results in the theory of random
graphs relate the local structure around a vertex in a random graph to
a branching process, see e.g.~\cite{Bol2,JLR}. More recently,
identification of the typical local limit structure for problems such
as random triangulations and quadrangulations and the random
assignment problem played an important part in establishing new and
exciting results, see e.g.~\cite{BS,AS2,AS}. 

In this paper we begin the study of typical local properties of a
model of random fixed 
width posets. This model was first studied by Brightwell and Goodall in \cite{BG} who derived an asymptotic formula for the number of such posets. Recall that an {\em antichain} of a poset $P$ is a subset
of mutually incomparable elements. The {\em width} of a poset
is defined to be the size of the largest antichain in it.

The width is a natural measure of the complexity of the set of partial
orders. Bounded width posets are of interest in
algorithms \cite{FT,DKMRV} and applications in artificial
intelligence \cite{SM,SWM} since
they can be specified in a compact manner \cite{DKMRV}.
For example, to specify a width-1 order, it is enough to give for each element
the minimal element larger than it. In fact, asymptotically, among the
set of all partial orders on $[n]=\{1, \cdots, n\}$, relatively few
have low width
\cite{BG,KR}.

\subsection{Results}
Let $\PP_{n,k}$ denote the set of labeled posets on vertex set
$[n]$ of width $k$. In \cite{BG}, the size of
$\PP_{n,k}$ is estimated asymptotically for fixed $k$ up to
factors polynomial in $n$. 

In this work, we analyze the {\em local structure} of random posets of
width $2$. In particular, we consider the asymptotic distribution of the
number of incomparable elements from the point of view of a random
element of the poset. As mentioned earlier, our work is best viewed in
light of recent success stories of the objective methods, or the use
of local weak limits, e.g.,~\cite{BS,AS2,AS} where identifying
limiting distributions of key properties shed light on problems that
were otherwise inaccessible. 

Our study is restricted to width $2$ posets. Even in this case the proofs are highly non-trivial while the case of posets of width $3$ 
 seems much harder both conceptually and technically. Another case where the analysis of width $2$ posets is possible and higher width posets seems impossible is~\cite{DT}.

To any poset $P$ on $[n]$ one can associate a graph $G_{P} = (V,E)$ with
vertices $V= [n]$ and edge set $E$ consisting of
those pairs $\{i, j\} \subset [n]$ which are \textit{not}
comparable. The graph $G_{P}$
is called the {\em incomparability graph} of $P$.
The incomparability graph $G_P$ captures much of the structure of
the poset $P$. For example, the width is the size of the maximum
clique of $G_P$ while the size of a maximum independent set is the
size of the largest chain, where a {\em chain} is a subset of
  elements of
the poset each pair of which is comparable.
Let $I_P(i)$ be the degree in $G_P$ of an element $i \in [n]$. Thus $I_P(i)$ is
exactly the number of elements incomparable to $i$ in $P$. We will
also refer to $I_P(i)$ as the {\em incomparability window} of $i$.

Let $\mu$ be the uniform distribution on $\PP_{n,2}$, let $P$ be
drawn according to $\mu$ and
let $U_n$ be a distribution that is uniform on
$[n]$ and independent of $\PP_{n,2}$ and $\mu$. We show that in the limit
as $n$ tends to infinity, $I_P(U_n)$ when
scaled appropriately converges in distribution to the height of a
Brownian excursion at a random time. We state the result more formally below.

Let $F_n$ be a sequence of cumulative distribution
functions corresponding to random variables $X_n$ and let
$F$ be the distribution function corresponding to a random variable
$X$. Recall (see for e.g. \cite{Dur}) that the
sequence $X_n$ is said to converge in distribution to
$X$ if for every $a \in
\mathbb R$ at which $F$ is continuous
$$    \lim_{n\rightarrow\infty} F_n(a) = F(a),$$
and denote it by
$$        X_n \stackrel{D}{ \Rightarrow} \ X.$$

Let $B^{\textrm{ex}}(t), 0 \le t \le 1$ denote the {\em Brownian excursion}
process on $[0,1]$ (see e.g. \cite{ItoMcKean}). Informally, a
Brownian excursion is obtained from the standard Brownian motion
$B(t)$ on $[0,1]$ by conditioning on the events $B(0) = B(1) = 0$ and
$B(t)>0$ for all $t \in (0,1)$. Let $U$ denote a uniform random
variable on $[0,1]$ independent of $B^{\textrm{ex}}(t)$.
Our main results are as follows.
\begin{theorem}\label{thm:distribution-of-window}
Let $P$ be a poset drawn from the distribution $\mu$. With the
notation defined above,
\begin{equation}
\frac{1}{\sqrt n} I_P(U_n) \overset{D}{\Rightarrow} \frac{1}{\sqrt
  2} B^{\textrm{ex}}(U)
\end{equation}
and follows the Rayleigh distribution with density $f(s)
= 8se^{-4s^2}, s>0$. 
\end{theorem}

\begin{remark}
It was observed by D. Aldous in \cite{Al} that $2B^{\textrm{ex}}(U)$ is
distributed as Rayleigh distribution with density $f(s) = se^{-s^2/2},
s> 0$.
\end{remark}

The {\em height} of a poset is the size of the longest chain in it. Let
$h(P)$ denote the height of the poset $P$.
\begin{theorem}\label{thm:height}
Let $P$ be a poset drawn from the uniform distribution $\mu$. Then
\[ \frac{h(P)  - n/2}{\sqrt n} \stackrel{D}{\Rightarrow} |Z|, \]
where $Z \sim N(0, 1/16).$
\end{theorem}

\subsection{Definitions and Notation}
Since we focus on fixed width posets with width $2$ on $n$ vertices,
to simplify notation 
henceforth we drop the subscripts and let $\PP$ denote the set
of width 2
posets on $[n]$. A poset $P = ([n], \prec)$ on ground set $[n] =
\{1,\cdots,n\}$ will be called a {\em labeled poset}. The set $\PP$ is
thus the set of labeled width-2 posets. Let $\bar P$ be the set of
posets {\em isomorphic} to $P$. Such a class $\bar P$ is called an {\em
  unlabeled poset}.

By Dilworth's theorem (see e.g. \cite{Bol1}), the ground set $[n]$ of a
width-2 poset $P= ([n], \prec)$ can be partitioned into two chains,
possibly in more than one way.
Let $P = (A \cup B, \prec)$ be a poset where $A \cup B$ is a
  partition of $[n]$ and $A$ and $B$ are chains in $P$.
Then $(A,B,\prec)$ is called a two-chain cover (of $P$), see for
example, Figure 1.

Let $C_1 = (A,B,\prec)$ and $C_2=(A',B',\prec)$ be two-chain covers (of
possibly different labeled partial orders).
Then $C_1$ and $C_2$ are isomorphic iff $|A| =
|A'|$ and there is a bijection taking each element of $A$ to
an element of $A'$
and each element of $B$ to an element
of $B'$ which preserves the relations of the partial order (for
example $C_1$ and $C_2$ in Figure 1). An isomorphism
class of two-chain covers consists of those two-chain covers which are
isomorphic to a specific two-chain cover. Such an
isomorphism class is called an {\em unlabeled two-chain cover} and the set of
all unlabeled two-chain covers of posets in $\PP$ will be denoted by
$\CC$. Note that isomorphic posets have the same sets of
unlabeled two-chain covers.

\vspace{0.2in}

\begin{pspicture}(-1.0,-2.7)(12,4)
\psset{dotsize=4pt 0}
\psset{linewidth=1pt}

\psline[showpoints=false,linecolor=black]{-}(0,.5 )(0,2)
\psline[showpoints=false,linecolor=black]{-}(1, 0)(1, 1)(1, 2)
\psline[showpoints=false,linecolor=black]{-}(0,.5 )(1,2)
\psline[showpoints=false,linecolor=black]{-}(0, 2 )(1,1)

\uput*{.2}[l](0,0.5){1}
\uput*{.2}[l](0,2){4}
\uput*{.2}[r](1,2){5}
\uput*{.2}[r](1,0){2}
\uput*{.2}[r](1,1){3}

\psdots[linecolor=black](1, 0)(1, 1)(1, 2)(0,.5 )(0,2)
\uput*{.2}[d](0.3,-1.5){$P = ([5],\prec)$}


\psline[showpoints=false,linecolor=black]{-}(3,2 )(3,3.5)
\psline[showpoints=false,linecolor=black]{-}(4, 1.5)(4, 2.5)(4, 3.5)
\psline[showpoints=false,linecolor=black]{-}(3,2 )(4,3.5)
\psline[showpoints=false,linecolor=black]{-}(3, 3.5 )(4,2.5)

\uput*{.2}[l](3,2){1}
\uput*{.2}[l](3,3.5){4}
\uput*{.2}[r](4,3.5){5}
\uput*{.2}[r](4,1.5){2}
\uput*{.2}[r](4,2.5){3}

\psdots[linecolor=red](3,2 )(3,3.5)
\psdots[linecolor=blue](4, 1.5)(4, 2.5)(4, 3.5)


\psline[showpoints=false,linecolor=black]{-}(6,2 )(6,3.5)
\psline[showpoints=false,linecolor=black]{-}(7, 1.5)(7, 2.5)(7, 3.5)
\psline[showpoints=false,linecolor=black]{-}(6,2 )(7,3.5)
\psline[showpoints=false,linecolor=black]{-}(6, 3.5 )(7,2.5)

\uput*{.2}[l](6,2){1}
\uput*{.2}[l](6,3.5){5}
\uput*{.2}[r](7,3.5){4}
\uput*{.2}[r](7,1.5){2}
\uput*{.2}[r](7,2.5){3}

\psdots[linecolor=red](6,2)(6,3.5)
\psdots[linecolor=blue](7, 1.5)(7, 2.5)(7, 3.5)

\psline[showpoints=false,linecolor=black]{-}(4.5, - .5 )(4.5,1)
\psline[showpoints=false,linecolor=black]{-}(5.5, -1)(5.5, 0)(5.5, 1)
\psline[showpoints=false,linecolor=black]{-}(4.5,-.5 )(5.5,1)
\psline[showpoints=false,linecolor=black]{-}(4.5, 1 )(5.5,0)

\uput*{.2}[l](4.5,-.5){4}
\uput*{.2}[l](4.5,1){1}
\uput*{.2}[r](5.5,1){5}
\uput*{.2}[r](5.5,-1){2}
\uput*{.2}[r](5.5,0){3}

\psdots[linecolor=red](4.5,-.5 )(4.5,1)
\psdots[linecolor=blue](5.5, -1)(5.5, 0)(5.5, 1)


\psline[showpoints=false,linecolor=black]{-}(10,.5 )(10,2)
\psline[showpoints=false,linecolor=black]{-}(9, 0)(9, 1)(9, 2)
\psline[showpoints=false,linecolor=black]{-}(10,.5 )(9,2)
\psline[showpoints=false,linecolor=black]{-}(10, 2 )(9,1)

\uput*{.2}[l](9,0){2}
\uput*{.2}[l](9,1){3}
\uput*{.2}[l](9, 2){5}
\uput*{.2}[r](10,.5){1}
\uput*{.2}[r](10,2){4}

\psdots[linecolor=red] (9, 0)(9, 1)(9,2)
\psdots[linecolor=blue](10,.5 )(10, 2)


\psline[showpoints=false,linecolor=black]{-}(13,.5 )(13,2)
\psline[showpoints=false,linecolor=black]{-}(12, 0)(12, 1)(12, 2)
\psline[showpoints=false,linecolor=black]{-}(13,.5 )(12,2)
\psline[showpoints=false,linecolor=black]{-}(13, 2 )(12,1)

\uput*{.2}[r](13,0.5){1}
\uput*{.2}[r](13,2){5}
\uput*{.2}[l](12,2){4}
\uput*{.2}[l](12,0){2}
\uput*{.2}[l](12,1){3}

\psdots[linecolor=red] (12, 0)(12, 1) (12, 2)
\psdots[linecolor=blue](13,.5 )(13,2)

\uput*{.2}[u](3.5,3.5){$C_1$}
\uput*{.2}[u](6.5,3.5){$C_2$}
\uput*{.4}[l](4.5, .5){$C_3$}
\uput*{.2}[u](9.5,2.1){$C'_1$}
\uput*{.2}[u](12.5,2.1){$C'_2$}
\uput*{.2}[l](2.5, -1.5){\Large $\bar{C}$}
\uput*{.2}[d](9.5, -1){\Large $\bar{C'}$}

\psdots[dotsize= 2pt] (4.7, -1.7) (5, -1.7)(5.3, -1.7)
\psdots[dotsize= 2pt] (11, -.7) (11.3, -.7)(11.6, -.7)

\psframe[linewidth=1pt,  linecolor=gray, framearc=.2,
](2.5, -2.0)(7.5,4.3)

\psframe[linewidth=1pt, linecolor=gray, framearc=.2,
](8.5, -1.0)(13.5,3)
\uput*{.2}[r](2,-2.5){$A$ consists of \textcolor{red}{red points}, $B$  consists of \textcolor{blue}{blue points} }

\end{pspicture}

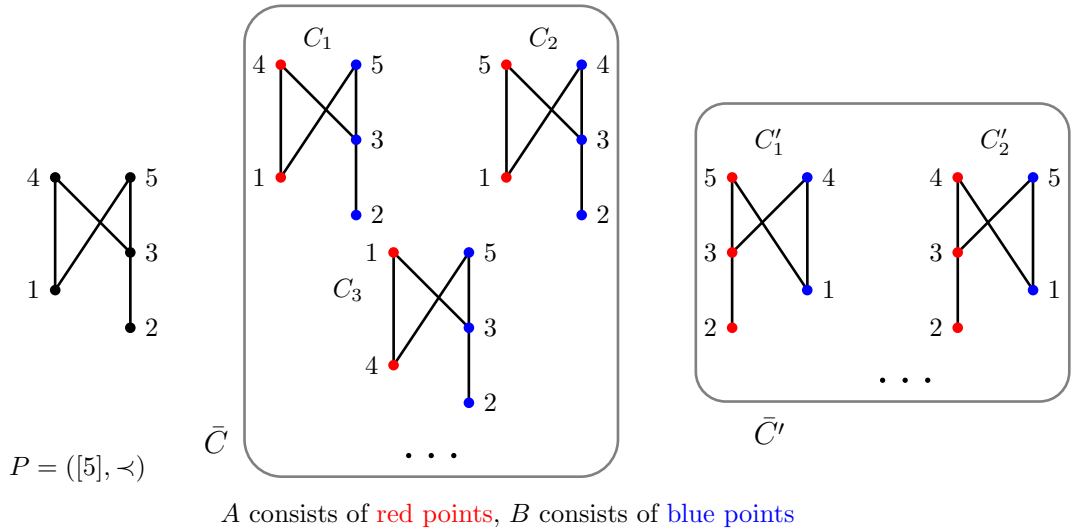
Figure 1: { \footnotesize In the figure, $P$ is an example of a width
 two labelled poset  on $5$ elements. $C_1, C_2, C_1', C_2'$ are all
 possible two-chain covers of $P$. Note  that $C_1$ and $C_2$ are
 isomorphic via the bijection $\vartheta :[5] \to [5]$ so that
 $\vartheta(i) = i$ for $1 \le i \le 3$ and $\vartheta(4) =
 5, \vartheta(5)=4$.  Thus they belong to the same unlabelled
 two-chain cover of $\bar C$ of $P$. 
 Similarly, $C_1'$ and $C_2'$ are also isomorphic while $C_1$ is not
 isomorphic to $C_1'$. Thus $C_1$ and $C_1'$ belong to two distinct
 two-chain covers $\bar C$ and $\bar C'$. Note that even though
 $C_3 \in \bar C$ is isomorphic to $C_1$ or $C_2$, it is a two-chain
 cover of some poset which is isomorphic to $P$ but not $P$ itself.
 }\\ 

The proof of Theorem \ref{thm:distribution-of-window} is based on the
construction of a bijection between
a set of unlabeled two-chain covers and constrained
pairs of non-hitting random walks on $\mathbb Z$. We
then reduce the analysis of the local
structure of the two-chain covers to the analysis of
the asymptotic
properties of these walks. This refines the approach of \cite{BG}
which used a bijection between $\mathcal{C}$ and a single constrained
one dimensional random walk to enumerate the size of $\PP$.

\subsection{Related Models of Random Posets}
The study of the asymptotic enumeration of posets was initiated
by Kleitman and Rothschild in \cite{KR}, where it was shown that
almost all posets have height 3 and width approximately
$\ffrac{n}{2}$.

There are a number of other models of random posets (see for
example the survey by Brightwell
\cite{Brightwell-survey}). For example, one can fix the fraction $0 < \alpha <
\frac{1}{2}$ of
comparable pairs and choose the poset uniformly from among those
posets with $\lfloor \alpha  n^2\rfloor$ comparable pairs
\cite{Dh}. The asymptotic analysis of the
evolution of the number of posets as a function of $\alpha$ is
given in \cite{KR2,PST}.
In a second model, one takes $d$ linear
orders on $[n]$ chosen uniformly at random and intersects them to
produce a poset
\cite{W,BB1}. The distribution of the heights of these posets was
analyzed in \cite{BW} and \cite{BDJ}.

\subsection{Open Problems}
In our opinion, an interesting open question is to obtain results
analogous to Theorem \ref{thm:distribution-of-window} for the
distribution of $I_P(i)$ for width-$k$ posets when $k>2$.
Unlike the width-2 case, it is not immediately clear how
the uniform measure on labeled width-$k$ posets gets
transferred to the geometric construction. Other open questions are
stated in Section \ref{sec:conclusions}.

\subsection{Organization of the Paper and Proof Structure}
The remainder of this paper is organized as follows. 
\begin{itemize} 
\item
In Section \ref{secn:reduce-to-one-factor} we show that the asymptotic
analysis for width-2 posets can be reduced to the analysis of width-2
posets with {\em one factor} (a factor of a poset will be defined later). 
The proof of this preliminary step utilizes some of the results of~\cite{BG}.
\item
In Section 3 we construct a bijection between two
chain covers with one factor and pairs of non-hitting random walks on
the lattice and relate statistics of the random walks to the sizes of
incomparability windows. This section contains the main contribution of the paper.
The mappings relating the different objects involved are non-trivial, and some of the proof require some novel 
combinatorial ideas, in particular Theorem~\ref{Lsym}. 

\item In Section 4, we give various estimates on the asymptotics of
the walks which allow
us to compute the limiting distribution of the incomparability
window. The challenge in the proof is that our main estimate, Theorem~\ref{thm:mapping-phi} gives implicit error bounds 
on the error in terms of the random objects we are trying to estimate. 

\item
Finally, given the machinery developed, it is straightforward to find the limiting distribution for the heights of
width-2 posets. This is done in Section 5.

\item
Section 6 lists a number of open problems.

\end{itemize}

\section{Reduction to Posets with One Factor}
\label{secn:reduce-to-one-factor}
In this section we reduce the problem of determining the limiting
distribution of the scaled incomparability window of a random element for a
random poset to the corresponding problem for posets with
exactly one factor
(see definition below). This reduction simplifies much of the
analysis in the subsequent sections.
\begin{definition}[\cite{BG}] Let $P \in \PP$ be a labeled width-2
  poset on $[n]$ and let
  $A_1,\cdots,A_k$ be the vertex sets of the connected components of the
  incomparability graph
  $G_P$. Further, let $P_i$ be the restriction of $P$ to
  $A_i$. The poset
  $P_i$ is called a {\em factor} of $P$.
\end{definition}
It can be verified that the factors of a poset can be ordered so that $P_1\prec
\cdots \prec P_k$, where $P_i \prec P_j$ means all elements of $P_i$
are less than all elements of $P_j$.
Let $\PP_1 \subseteq \PP$ be the set of posets with one factor and let
$\BbbP_\texttt{pos}$ be the uniform distribution over $\PP_1$.
\begin{lemma}
\label{lem:reduction-to-one-factor}
Let $Q$ be a poset chosen uniformly at random according to
  {\em $\BbbP_\texttt{pos}$} and let $U_n$ be uniform over the ground set $[n]$
  independent of $Q$. Let $B^{\text{\emph{ex}}}$ be the standard Brownian
  Excursion on
  $[0,1]$ and $U \sim$ Uniform $[0,1]$ independent of
  $B^{\text{\emph{ex}}}$. Suppose that
\begin{equation}
 \label{eq:convergence-for-one-factor}
 \frac{I_{Q}(U_n)}{\sqrt n}  \stackrel{D}{ \Rightarrow} \frac{1}{\sqrt
 2} B^{\text{\emph{ex}}}(U).
 \end{equation}
Then, if $P$ is a poset chosen according to $\mu$ and $U'_n$
is uniform over the ground set $[n]$ and independent of $P$
\begin{equation}
 \label{eq:convergence-for-general case}
 \frac{I_{P}(U'_n)}{\sqrt n}  \stackrel{D}{ \Rightarrow} \frac{1}{\sqrt
 2} B^{\text{\emph{ex}}}(U).
 \end{equation}
\end{lemma}
Theorem
\ref{thm:distribution-of-window} is then
reduced to proving the hypothesis
(\ref{eq:convergence-for-one-factor}), and will be proved in
Section \ref{sec:convergence-for-one-factor}.
The proof of Lemma \ref{lem:reduction-to-one-factor} relies on the
following result of Brightwell and Goodall
(see page 329 of \cite{BG}) showing that asymptotically, most
width-2 posets have one large factor.
\begin{lemma}[\cite{BG}]
\label{lem:BG-large-factor}
Let $P$ be chosen uniformly at random from $\PP$. Then
\begin{equation}
\mathbb{P}(G_P \ \text{has \ a \ factor of size }  \geq n- \log n)
\rightarrow 1
\end{equation}
as $n \rightarrow \infty$.
\end{lemma}
\begin{proof}[Proof of Lemma \ref{lem:reduction-to-one-factor}]
Let $U_n$ be chosen uniformly from the ground set $[n]$.
Consider the random variable $I_P(U_n)$ and let $M$ denote the
largest connected component of $G_P$. By Lemma
\ref{lem:BG-large-factor}, for $0<a<b$,
\begin{equation}
\mathbb P\left(  I_P(U_n) \in (a\sqrt{n}, b\sqrt{n}]\right) = \mathbb
P\left( I_P(U_n) \in (a\sqrt{n}, b\sqrt{n}] \cap \{|M| \geq n -
  \log n\}\right) + o(1).
\end{equation}
Conditioning on the size of $M $, we have that $\mathbb
P\left(U_n \in M \big | |M|\right) = n^{-1}|M|$, so
we may further write
\begin{align*}
\mathbb P\left(I_P(U_n) \in (a\sqrt{n}, b\sqrt{n}]\right) =
\mathbb P \left( I_P(U_n) \in (a\sqrt{n}, b\sqrt{n}] \big | U_n
  \in M, |M|
  \geq n - \log n \right) + o(1).
\end{align*}
For each $m \geq n- \log n$, let $F_m=\{ U_n \in M \} \cap
\{|M|=m\}$. Since with high probability $U_n$ is chosen from
$M$ when $|M| \geq n -\log n$,
\begin{align*}
\mathbb P\left(  I_P(U_n) \in (a\sqrt{n}, b\sqrt{n}]\right) =
\sum_{m\geq n- \log n} \mathbb P\left(  I_P(U_n) \in (a\sqrt{n},
  b\sqrt{n}]| F_m\right)\mathbb P(F_m) + o(1).
\end{align*}

By construction, the distribution of $I_P(U_n)$ conditioned on $F_m$
is the same as the distribution of $I_{Q}(U_m)$ where $Q$ is a
random poset with one factor on the ground set $[m]$ and $U_m$
is uniform on $[m]$. In other words,
\begin{align*}
\mathbb P\left(  I_P(U_n) \in (a\sqrt{m}, b\sqrt{m}]\right | F_m) =
\mathbb P\left(  I_Q(U_m) \in (a\sqrt{m}, b\sqrt{m}]\right).
\end{align*}

Thus from
\eqref{eq:convergence-for-one-factor},
for all $\epsilon>0$, there is $m(\epsilon)$ so that for
$m>m(\epsilon)$,
\begin{equation}
\left|\mathbb P\left(  I_P(U_n) \in (a\sqrt{m}, b\sqrt{m}]|
    F_m\right)-\mathbb P\left( \frac{1}{\sqrt 2} B^{\text{ex}}(U) \in
    (a, b]\right) \right| \leq \epsilon.
\end{equation}

Note that this holds for any pair $a, b$ (not just the particular pair of interest) and that the function $P\left( \frac{1}{\sqrt 2} B^{\text{ex}}(U) \in (a, b]\right)$ is (Lipschitz) continuous in $a, b$.
Thus, choosing $n$ large enough so that $n - \log n \geq m(\epsilon)$ and
$o(1) \leq \epsilon$ we have
\begin{align*}
\left|\mathbb P\left(  I_P(U_n) \in (a\sqrt{n}, b\sqrt{n}]\right)-
  \mathbb P\left( \frac{1}{\sqrt 2} B^{\text{ex}}(U) \in (a,
    b]\right)\right| & \ \leq
\epsilon \sum_{m\geq n- \log n}\mathbb P(F_m) + o(1) \\
& \ \leq 2 \epsilon
\end{align*}
and the claim follows.
\end{proof}

As a consequence of Lemma \ref{lem:reduction-to-one-factor},
we can henceforth restrict ourselves to working over $\PP_1$, the set of
labeled posets with just a single factor.

\section{Combinatorial Mappings} \label{sec:combinatorial-maps}
In this section we construct a correspondence between posets $\PP_1$
and a random walk representation useful for our probabilistic analysis. We
then relate certain statistics of the random walks to the size
of the incomparability window. The
idea behind the correspondence is essentially due to Brightwell and
Goodall, but requires some modification in order to describe the
incomparability window.

Denote the set of unlabeled two-chain covers with
one factor by $\CC_1$. We say that $\bar C \in \CC_1$ is an
{\em unlabeled two-chain cover of a poset $P$} if there is $C' \in
\bar C$ which is a two-chain cover of $P$.

For each unlabeled two-chain cover $\bar C$, let $C = (A,B,\prec)$ be the
element of $\bar C$
with $A = \{1, \cdots,k\}$ and $B = \{k+1, \cdots,n\}$,
where $k$ is the size of the first chain, so that the chains
$A$ and $B$ are ordered in the 
natural way. The
two-chain cover $C$ will be the canonical representative of $\bar C$ and
we will identify elements of $\CC_1$ with their canonical
representatives.  By abuse of notation, by $C = (A,B,\prec) \in \CC_1$ we
will mean the element of $\CC_1$ which is the isomorphism class of
$C$.

The following correspondence between $\PP_1$ and $\CC_1$
was given in \cite{BG}. Let $P \in \PP_1$ and define
\begin{align}\label{eq:psi}
\Psi(P) = \{C \in \CC_1 \ : \ C \text{ is an unlabeled two-chain
  cover of } P   \}.
\end{align}
Note that $|\Psi(P)| \leq 2$.
It can be verified that $\Psi$
associates disjoint sets of covers with non-isomorphic posets and
associates the
same set of covers with any two isomorphic posets.
For any unlabeled two-chain cover $C$, define $\Psi^{-1}(C)$ to be the
set of posets $P$ in $\PP_1$ such that $C$ is an unlabeled two-chain
cover of $P$.
The random walk representation we study is as follows.
\begin{definition}\label{def:non-hitting-walks}
A pair of {\em non-hitting walks} $(V, W)$ of length $n$ on the
integer lattice is defined to be a pair of walks where
\vspace{-3.6pt}
\begin{enumerate}
\item $V(0) = W(0) = 0.$
\item For every $0 \leq t < n$, $|V(t+1) - V(t)| =|W(t+1) - W(t)|= 1.$
\item $V(n) = W(n)$.
\item For every $1 \leq t < n$, $V(t) > W(t)$.
\end{enumerate}
\end{definition}
Let $\BB_n$ be the set of non-hitting pairs of walks of length
$n$. For $(V, W) \in \BB_n$, define
the height function $H_{V,W}:[n]
\rightarrow \{0, 2,4,\cdots\}$ by \[ H_{V,W}(k) := V(k) - W(k) \ \ \ \
\ \ 1 \leq k \leq n.\]
\begin{definition}
Let $P = ([n],\prec)$ be a labeled poset. Let us define the function $\tau_P : [n]
\rightarrow [n]$ by
\begin{align}\label{eq:tau}
\tau_P(x) = |\{y \in [n] \ : \ y\preceq x\}|.
\end{align}
\end{definition}
The first result of this section relates the incomparability window
$I_P(x)$ for an element $x$ of a poset $P \in \PP_1$ to
$H_{V,W}(\tau_P(x))$ for walks $(V,W) \in \BB_n$ associated with $P$.

The set of width-2 posets and non-hitting
walks are associated as follows.
We will construct a bijection $\Gamma$ between $\CC_1$ and the set
of non-hitting walks $\BB_n$ in Section \ref{S:Bij}.
Given $\Gamma$, let
\begin{align}\label{eq:phi}
\Phi(P) = \Gamma(\Psi(P))
\end{align}
be the image of $\Psi(P)$ under $\Gamma$.
The incomparability window can then be related to
the height function $H_{V,W}$ as follows.
\begin{theorem}\label{thm:mapping-phi}
Let $P \in \PP_1$, $(V,W) \in \Phi(P)$ and let $x \in [n]$ be an
  element of the ground set. Then,
\begin{align*}
|I_{P}(x)-H_{V,W}(\tau_P(x))| \leq Err(P, x)
\end{align*}
where
\begin{align*}
Err(P,x) : =
 |V(\tau_P(x) + I_P(x)) - V(\tau_P(x))| + |W(\tau_P(x) + I_P(x)) -
 W(\tau_P(x))|.
\end{align*}
\end{theorem}

At a high level, the argument for why the above bound is useful is as
follows. We will show that under $\Phi$ the uniform distribution on
$\PP_1$ induces the uniform distribution over
$\BB_n$. For $(V,W)$ chosen uniformly from
$\BB_n$, we expect that under diffusive scaling,  $H_{V,W}$
converges to a Brownian excursion.
Our goal is to show that for a uniformly random $P$ and a uniformly
random element $U_n$, the
scaled error $Err(P,U_n)/ \sqrt n$ goes to 0 in probability.
In Section \ref{sec:convergence-for-one-factor} we show that with high
probability the size of
the incomparability window is bounded by $n^{2/3}$. Therefore, bounding
the scaled error reduces to bounding
the maximum fluctuation of each the walks $V$ and $W$ in
$n^{2/3}$ steps.
We will show that the
differences $|V(\tau_P(U_n) + I_P(U_n)) - V(\tau_P(U_n))|$ are at most
$n^{1/3}$ with high probability and hence the scaled error goes
to 0.

Define the joint distribution $\BbbP$ on $\PP_1 \times \CC_1 \times
\BB_n$ as follows. Let $P \sim \BbbP_\texttt{pos}$ be chosen uniformly from
$\PP_1$. Given $P$, choose $C$ uniformly from $\Psi(P)$. Given $P$ and
$C$, let $(V,W) = \Gamma(C)$
with probability $1$. Let $\BbbP_\texttt{walk}$ and $\BbbP_\texttt{cov}$ be the marginals of
$C$ and $(V,W)$, i.e. for every $C \in \CC_1$
\begin{align}
\BbbP_\texttt{walk}(C) = \sum_{P: C \in \Psi(P)}
\frac{\BbbP_\texttt{pos}(P)}{|\Psi(P)|}
\end{align}
and $(V,W) \in \BB_n$
\begin{align}
\BbbP_\texttt{cov}(V,W) = \sum_{P: (V,W) \in \Phi(P)}
\frac{\BbbP_\texttt{pos}(P)}{|\Phi(P)|}.
\end{align}

\begin{lemma}[\cite{BG}]\label{lem:mu0-uniform-over-T0}
Let  {\em$\BbbP_\texttt{walk}$} be defined as above. The distribution {\em $\BbbP_\texttt{walk}$} is uniform
over $\CC_1$.
\end{lemma}
\begin{lemma} \label{lem:sig}
Let  {\em $\BbbP_\texttt{cov}$} be defined as above. Then {\em $\BbbP_\texttt{cov}$} is the
uniform
distribution over $\BB_n$.
\end{lemma}

The proof is postponed until after the construction of the bijection
$\Gamma$ (Proposition \ref{prop:gamma-is-bijection}).

\begin{definition}
Let $(V,W)  = \Gamma(C)$. For $x \in
[n]$, define the function $\tau_{V,W}(x) = |\{y \ : \ y \preceq_{C}
x\}|$  where $y \preceq_C x $ denotes that $y \preceq x$ in the partial
order defined by $C$.
\end{definition}

\begin{lemma}\label{lem:tauP-tauT}Let $(P,(V,W))$ be chosen according to the
  two-dimensional marginal of $\BbbP$ on $\PP_1 \times
  \BB_n$. Let $U_n$ be a uniform element of the ground set
  $[n]$ independent of $(P,(V,W))$. Then, there exists a random
  variable $U_n'$ which is uniform over
  $[n]$ and independent of $(V,W)$ such that
\begin{align}
H_{V,W}(\tau_P(U_n)) \stackrel{d}{=} H_{V,W}(\tau_{V,W}(U_n'))
\end{align}
\end{lemma}

\begin{proof}
Since $P$ and $(V, W)$ are chosen according to the two-dimensional
marginal of $\BbbP$, $(V, W) \in \Phi( P )$. Let $C= \Gamma^{-1}((V,
W))$.  In particular $C \in \Psi(P)$ where we are using the convention
that $C$ is the canonical
representative of an unlabeled two-chain cover $\bar C$ of $P$.   Since
$C \in \Psi(P)$,  for any (labeled) two-chain cover  $C_P \in \bar C$ of
$P$ there is a (unique) two-chain isomorphism $\ga$ between  $C_P$ and
$C$. Consequently,
\begin{align*} \tau_P(x) &= |\{ y \in [n] : y \preceq x \} | \\
&= |\{ y \in [n] : y \preceq_{C_P} x \} | \\
&=   |\{ y \in [n] : y \preceq_{C} \ga(x) \} |  = \tau_{V, W}(\ga(x)).
\end{align*}
Given $P$ and $(V, W) \in \Phi( P )$. Let
\[ U_n' = \ga(U_n).\]
Note that though the map $\alpha$ depends on $P, (V, W)$ and the
choice of $C_P$, since $\mathcal U_n$ is uniform over $[n]$ and
independent of everything else, $U_n'$ is also uniform over $[n]$ and
independent of $P$ and $(V, W)$. The assertion of the lemma follows.
\end{proof}

In order to formalize the high level argument above we would like to relate the
distribution of $H_{V,W}(\tau_{V,W}(U_n))$ for a uniformly
random element $U_n$ to the distribution of $H_{V,W}$ evaluated
at a
uniformly random time. This is the main result of the section.
\begin{theorem}[Symmetrization]\label{Lsym} Let $(V,W)$ be chosen at
  random from $\BB_n$ according to {\em $\BbbP_\texttt{cov}$}, and let $C =
  \Gamma^{-1}(V,W) = (A,B,\prec)$.
Let $U_n$ be a
uniform element from the ground set $[n]$ and let
$U_n'$ be uniform on
$ \{1,2, \ldots, n\}$, both independent of $(V,W)$.  Then for every $m \in
2\mathbb{Z}_{\geq}$,
\[ \mathbb P(H_{V,W}(\tau_{V,W}(U_n)) =m) = \mathbb
P ( H_{V,W}(U_n') = m). \]
\end{theorem}

\subsection{The Bijection $\Gamma$ from Two-chain Covers to Non-hitting Walks}
\label{S:Bij}
Next, we construct a bijection $\Gamma$
from $\CC_1$ to the set of
pairs of walks $\BB_n$.
Our construction deviates from that of \cite{BG}
since we construct a bijection to pairs
of walks rather than a single walk.

Given a two-chain cover $C = (A,B,\prec)$ in $\CC_1$, let
$(\lambda, \delta)$ denote the pair of
total orders,
defined by adding to the poset the following
relations: if $a \in
A, b \in B$ are not comparable in $C$, then we set $a\prec b$ in $\lambda$ and set
$a\succ b$ in $\delta$. Following the terminology of \cite{BG}, we call
  $(\lambda,\delta)$ the {\em
  greedy pair} of orders associated with $C$. The relations of the two
chain cover can be reconstructed by setting $a \prec  b \ (\mathrm{resp.}
\succ)$ if $a$ precedes (resp. follows) $b$ in both of the greedy total
orders.

For $c \in [n]$, denote by $\texttt{rank}_\lambda(c)$ the number of elements
less than or equal to $c$ in $\lambda$ and define $\texttt{rank}_\delta(c)$
similarly. For $1 \leq i \leq n$, let $\lambda_i$ (resp. $\delta_i$)
denote the element in position $i$ in $\lambda$ (resp. $\delta$).

\begin{proposition}[\cite{BG}]
  \label{prop:characterization-of-greedy-orders} A pair of orders
  $(\lambda,\delta)$ on
  $[n]$ is a greedy pair for some two-chain cover $C = (A,B,\prec)$ iff for
  every $a \in A$ and
  $b \in B$, if {\em $\texttt{rank}_\delta(a) < \texttt{rank}_\delta(b)$} then {\em $\texttt{rank}_\lambda(a) < \texttt{rank}_\lambda(b)$}.
\end{proposition}

Given the pair of greedy orders $(\lambda, \delta)$ obtained from $C$,
we define a pair of
walks $(V, W)$ as follows:
Start by setting $V(0) = W(0) = 0$. In the $i$-th step, the walk $V$
(resp. $W$) takes a step up if the $i$-th element of the total order
$\lambda$ (resp. $\delta$) belongs to $A$, and otherwise it
takes a step down (see for example Figure 2).  In other words

\begin{align*}
V(i) = \left\{ \begin{array}{ll}
V({i-1})+1 & \text{ if}
 \ \lambda_i \in A \\
V({i-1})-1 & \text{ if} \ \lambda_i \in B
\end{array} \right.
\hspace{0.25in}
W(i) = \left\{ \begin{array}{ll}
W({i-1})+1 & \text{ if}
\ \delta_i \in A\\
W({i-1})-1 & \text{ if} \ \delta_i \in B.
\end{array} \right.
\end{align*}

For $1\le k \le n$, let $\pi_k(V) := V(k) - V(k-1)
\in \{-1, 1\}$ denote the $k$-th step of $V$ and define $\pi_k(W)$ analogously.

\begin{minipage}[b]{0.3\linewidth} 
\begin{pspicture}(-1.0,-0.5)(5,5)
\pscurve[showpoints=true,linecolor=black]{-}(0,1)(0,2) (0,3)(0,4)
\pscurve[showpoints=true,,linecolor=black]{-}(1.2,1.5)(1.2,2.5)
(1.2,3.5)

\uput*{.2}[l](0,1){$a_1$}
\uput*{.5}[d](0,1){ $A$}
\uput*{1}[d](1.2,1.5){ $B$}
\uput*{.2}[l](0,2){$a_2$}
\uput*{.2}[180](0,3){$a_3$}
\uput*{.2}[180](0,4){$a_4$}
\uput*{.2}[r](1.2,1.5){$b_1$}
\uput*{.2}[0](1.2,2.5){$b_2$}
\uput*{.2}[0](1.2,3.5){$b_3$}

\psline[linecolor=black]{-}(0,3)(1.2,1.5)
\psline[linecolor=black]{-}(0,4)(1.2,2.5)
\psline[linecolor=black]{-}(0,2)(1.2,3.5)

\pscurve[showpoints=true,linecolor=black]{-}(3,-1)(3,0)(3,1 )(3,2 ) (3,3)(3,4)(3,5)
\pscurve[showpoints=true,linecolor=black]{-}(4,-1)(4,0)(4,1 )(4,2 ) (4,3)(4,4)(4,5)

\uput*{.2}[l](3,-1){$a_1$}
\uput*{.2}[l](3,0){$a_2$}
\uput*{.2}[l](3,1){$b_1$}
\uput*{.2}[l](3,2){$a_3$}
\uput*{.2}[l](3,3){$b_2$}
\uput*{.2}[l](3,4){$a_4$}
\uput*{.2}[l](3,5){$b_3$}
\uput*{.5}[d](3,-1){$\lambda$}

\uput*{.2}[r](4,-1){$b_1$}
\uput*{.2}[r](4,0){$b_2$}
\uput*{.2}[r](4,1){$a_1$}
\uput*{.2}[r](4,2){$a_2$}
\uput*{.2}[r](4,3){$b_3$}
\uput*{.2}[r](4,4){$a_3$}
\uput*{.2}[r](4,5){$a_4$}
\uput*{.5}[d](4,-1){$\delta$}

\end{pspicture}
\end{minipage}
\hspace{0.5cm} 
\begin{minipage}[h]{0.7\linewidth}
\begin{pspicture}(-4.0,-2.0)(4,8)
\psaxes[linewidth=.5pt,labels=none,ticks=none]{<->}(-1,3.5)(-1.4,1)(5,6)
\uput*{.2}[230](-1,3.5){O}
\psline[showpoints=true]{-}(-1,3.5)(-.25,4.25)
\uput*{.2}[u](-.25,4.25){a$_{1}$}
\psline[showpoints=true]{-}(-.25,4.25)(.5,5)
\uput*{.2}[u](.5,5){a$_{2}$}
\psline[showpoints=true]{-}(.5,5)(1.25, 4.25)
\psline[showpoints=true]{-}(1.25, 4.25)(2, 5)
\uput*{.2}[u](2,5){a$_{3}$}

\psline[showpoints=true]{-}(2, 5)(2.75, 4.25)
\psline[showpoints=true]{-}(2.75, 4.25)(3.5, 5)
\uput*{.2}[u](3.5,5){a$_{4}$}
\psline[showpoints=true]{-}(3.5, 5)(4.25, 4.25)

\uput*{.2}[u](.2,5.5){V}
\uput*{.2}[d](.2, 1.2){W}

\psline[showpoints=true]{-}(-1,3.5)(-.25,2.75)
\uput*{.2}[d](-.25,2.75){b$_{1}$}
\psline[showpoints=true]{-}(-.25,2.75)(.5, 2)
\uput*{.2}[d](.5,2){b$_{2}$}
\psline[showpoints=true]{-}(.5, 2) (1.25, 2.75)
\psline[showpoints=true]{-}(1.25, 2.75)(2, 3.5)
\uput*{.2}[135](2,3.5){a$_{2}$}
\psline[showpoints=true]{-}(2, 3.5)(2.75,2.75)
\uput*{.2}[d](2.75,2.75){b$_{3}$}
\psline[showpoints=true]{-}(2.75,2.75)(3.5, 3.5)
\uput*{.2}[320](3.5,3.5){a$_{3}$}
\uput*{.2}[d](5,3.5){$t$}
\psline[showpoints=true]{-}(3.5, 3.5)(4.25, 4.25)

\psline[showpoints=true, linecolor= cyan,linestyle =dashed]{<->}(1.9, 4.95)(3.4, 3.5)
\psline[showpoints=true, linecolor= magenta, linestyle=solid]{|<->|}(2, 2)(3.5, 2)
\uput*{.2}[d](2.75,2){\scriptsize$I_C(a_3) $}
\psline[showpoints=false, linecolor= yellow]{->}(2, 1)(2, 3.5)
\psline[showpoints=false, linecolor= yellow]{->}(3.5, 1)(3.5, 3.5)
\uput*{.2}[d](2,1){\scriptsize $t_V(3)$}
\uput*{.2}[d](3.5,1){\scriptsize $t_W(3)$}


\end{pspicture}
\end{minipage}
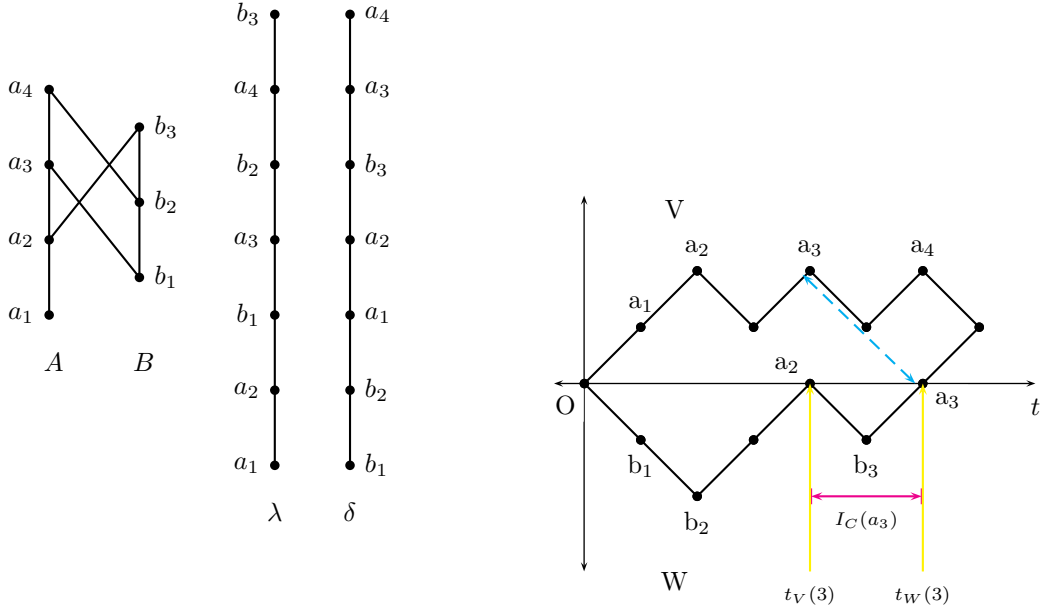
Figure 2: {\footnotesize (Left) A two-chain covering $(A,B, \prec)$ of a poset
with $7$ elements. (Middle) The corresponding  left and right greedy orders $(\lambda, \delta)$. (Right) The associated pair of non-crossing walks
$(V,W)$. The incomparability window for the element $a_3$ is
illustrated.}

\begin{proposition}\label{prop:map-to-non-hitting-walks}
Let $C \in \CC_1$. Then $\Gamma(C) = (V,W)$ is a pair of
non-hitting walks.
\end{proposition}
\begin{proof} Recall that non-hitting walks satisfy four properties (see
  Definition \ref{def:non-hitting-walks}). Properties $(1), \ (2)$ and
  $(3)$ follow immediately by the definition of $\Gamma$.
Suppose that $(4)$ does not hold, so that in particular, there is a
smallest $k$,  $0
< k < n$ such that $V(k) = W(k)$. The claim is that then $C$ must have
more than one factor, a contradiction.

By the definition of $\lambda$ and $\delta$, every element $a \in A$
satisfies $\texttt{rank}_\lambda(a) \leq \texttt{rank}_\delta(a)$ and conversely, every
element of $b \in B$ satisfies $\texttt{rank}_\delta(b) \leq \texttt{rank}_\lambda(b)$. Since $V(k)
= W(k)$, by construction
the sets $S_\lambda = \{\lambda_1, \cdots,\lambda_k\}$ and
$S_\delta = \{\delta_1,\cdots,\delta_k\}$ consist of the same numbers
of elements
from $A$ and $B$. Now consider any $\lambda_i \in B \cap S_\lambda$ so
that $\lambda_i=b$, i.e., $\texttt{rank}_\lambda(b)=i$. We have that $\texttt{rank}_\delta(b) \leq
\texttt{rank}_\lambda(b)=i \leq k$, and therefore, $b \in S_\delta$. A
symmetric argument for $\delta_j \in A \cap S_\delta$ implies
$\delta_j \in S_\lambda$, and hence, $S_\lambda = S_\delta  = S$.

Since $|S| = k < n$, there exists an element $x \notin S$ appearing in both
$\lambda$ and $\delta$ above the elements of $S$. The element
$x$ must therefore be in a different component of the incomparablity
graph of the poset $P$ defined by $C$ than the
elements of $S$, implying that $C$ has more than one factor.
\end{proof}
\begin{proposition}\label{prop:gamma-is-bijection} The map $\Gamma: \CC_1
  \rightarrow
  \BB_n$ is
  a bijection.
\end{proposition}
\begin{proof}
By the definition of the map $\Gamma$, distinct two-chain covers are
mapped to non-isomorphic greedy orders. Moreover, non-isomorphic
greedy orders are
mapped to distinct walks $(V, W)$ and by Proposition
\ref{prop:map-to-non-hitting-walks} the walks
are non-hitting. Hence $\Gamma$ is
one-to-one. It remains to show that any pair of walks corresponds to
some two-chain cover.

From a pair of walks $V,W \in \BB_n$ define a pair of total orders
$(\lambda,\delta)$ on $[n]$ and sets $A$ and $B$, where $A \cup B =
[n]$ as follows. Let $A = \{1,
\cdots, k\}$ and $B = \{k+1, \cdots,n\}$, where $k
= |\{t \ : \ \pi_t(V) = +1\}| = |\{t \ : \ \pi_t(W) = +1\}|$ and the equality follows since $V(n) = W(n)$.
The greedy orders are defined as follows.
Let $V_{+ } = (i_1, \cdots,i_{|A|}) $ be the (increasing) subsequence of indices
so that $\pi_{i_{\ell}}(V) = +1$.  Analogously, let us define the
subsequences of indices $V_- = (i_{|A|+1},\cdots,i_n)$, $W_+=(
j_{1},\cdots,j_{|A|})$, and $W_- =( j_{|A|+1},\cdots,j_n)$.
The linear orders $\lambda$ and $\delta$ are then defined by the rankings
$\lambda_{i_\ell} = \ell$ and $\delta_{j_\ell} = \ell$ for each $\ell
\in [n]$.

By property $(4)$ of the walks $(V,W)$, there are no elements $a \in
A$ and $b \in B$ such that $\texttt{rank}_\delta(a) < \texttt{rank}_\delta(b)$ but $\texttt{rank}_\lambda(a) >
\texttt{rank}_\lambda(b)$. Therefore the pair $(\lambda,\delta)$ satisfies the
condition of Proposition
\ref{prop:characterization-of-greedy-orders} and is a greedy pair of
linear orders
corresponding to some two-chain cover.
Hence, the map $\Gamma$ is a bijection.
\end{proof}
\vspace{-4pt}
Lemma \ref{lem:sig} is now an easy corollary.
\vspace{-4pt}
\begin{proof} [ Proof of Lemma \ref{lem:sig} ]
Let $P$ be a
  width-2 poset in $\PP_1$. Let $\Psi$ and $\Phi$ as defined in
  \eqref{eq:psi} and \eqref{eq:phi}. Suppose that $(V_1,W_1)$ and
  $(V_2,W_2)$ are walks in
  $\BB_n$. Let $C_1 = \Gamma^{-1}(V_1,W_1)$  and $C_2 =
  \Gamma^{-1}(V_2,W_2)$. Then,
\begin{align*}
\BbbP_\texttt{cov}(V_1,W_1) =  \sum_{P \in \Phi^{-1}(V_1,W_1)}
\frac{\BbbP_\texttt{pos}(P)}{|\Phi(P)|} = \sum_{P \in \Psi^{-1}(C_1)}
\frac{\BbbP_\texttt{pos}(P)}{|\Psi(P)|} \\ = \sum_{P \in \Psi^{-1}(C_2)}
\frac{\BbbP_\texttt{pos}(P)}{|\Psi(P)|} = \sum_{P \in
  \Phi^{-1}(V_2,W_2)}\frac{\BbbP_\texttt{pos}(P)}{|\Phi(P)|}   = \
\BbbP_\texttt{cov}(V_2,W_2)
\end{align*}
The third equality follows from Lemma \ref{lem:mu0-uniform-over-T0}
and the fourth by Proposition \ref{prop:gamma-is-bijection}.
\end{proof}

\begin{definition}
Let $C = (A,B, \prec) \in \CC_1$. For $i \in [n]$, define the
incomparability window of $i$ in $C$, $I_C(i):=I_P(i)$ where $P$ is
the poset defined by $C$.
\end{definition}

The elements of $C=(A, B,\prec) \in \CC_1$ can be  labeled $a_1\prec
\cdots\prec  a_k$ and $b_1\prec \cdots \prec b_{n-k}$.
For each $1 \leq i \leq k$, let $t_V(i)$ be the first time at
which the walk $V$ has taken $i$
steps up and let $t_W(i)$ be the first time $W$
has taken $i$ steps up. Similarly, for each $1 \leq j \leq n-k$ let
$s_V(j)$ and $s_W(j)$ be the first times
that $V$ and $W$ have taken $j$ downward steps respectively.

\begin{lemma}\label{lem:window-is-intercept}
Let $n$ be fixed and let $C=(A,B,\prec)$ be a two-chain cover in $\mathcal
C_1$.  Then we
have the following identifications:
\begin{enumerate}
\item The size of the incomparability window for $a_i$ satisfies
\[ I_C(a_i) = t_W(i)-t_V(i).\]
Similarly,  for an element $b_j$,
\[ I_C(b_j) = s_V(j) - s_W(j).\]
\item
Alternatively,
\begin{equation}
\label{eq:size}
I_C(a_i) = H_{V,W}(t_V(i))  + W(t_V(i)) - W(t_V(i) + I_{T}(a_i))
 \end{equation}
 and
 \begin{equation}
   \label{eq:size2}
   I_C(b_j) = V(s_W(j)+I_T(b_j)) - V(s_W(j)) + H_{V,W}(s_W(j)).
 \end{equation}
 \end{enumerate}
\end{lemma}

\begin{proof}
By symmetry we may restrict attention to the statements for $a_i \in
A.$  We first prove part $(1)$.
At time $t_V(i)$ when $V$ has taken $i$ steps `up', the set
$\{\lambda_1,\cdots,\lambda_{t_V(i)}\}$ consists of $a_i$ and
all the elements below $a_i$. Similarly, when $W$ has taken $i$ steps
`up', at

\begin{minipage}[b]{0.3\linewidth}
\begin{pspicture}(-1.0,-1.0)(2,4)
\psset{dotsize=5pt 0}
\psset{linewidth=1.2pt}

\pscurve[showpoints=true,linecolor=black]{-}(0,1 )(0,3 ) (0,5)
\pscurve[showpoints=true,linecolor=black]{-}(2, .5)(2, 2.5) (2, 4.5)

\uput*{.2}[l](0,1){$a_1$}
\uput*{1}[d](0,1){ $A$}
\uput*{.5}[d](2,.5){ $B$}
\uput*{.2}[l](0,3){$a_2$}
\uput*{.2}[180](0,5){$a_3$}

\uput*{.2}[r](2,.5){$b_1$}
\uput*{.2}[0](2,2.5){$b_2$}
\uput*{.2}[0](2,4.5){$b_3$}

\psline[linecolor=black]{-}(0,3)(2,4.5)
\psline[linecolor=black]{-}(0,5)(2, .5)
\psline[linecolor=blue, linestyle =dashed]{-}(0,3)(2,2.5)
\psline[linecolor=blue, linestyle =dashed]{-}(0,3)(2,.5)




\end{pspicture}
\end{minipage}
\hspace{1.5cm} 
\begin{minipage}[b]{0.7\linewidth}
\begin{pspicture}(-1.0,-4)(12,4)
\psset{dotsize=5pt 0}
\psset{linewidth=1.2pt}
\psgrid[gridwidth=.5pt,subgriddiv=1,subgridwidth=.1pt,gridcolor=black,gridlabels=0](0,-3)(7, 3)
\psline[showpoints=false, linecolor =red]{-}(0,0)(1,1)

\uput*{.3}[150](1,1){a$_{1}$}
\psline[showpoints=false, linecolor =red]{-}(1,1)(2,2)
\uput*{.3}[150](2,2){a$_{2}$}
\psline[showpoints=false,  linecolor =red]{-}(2,2)(3,1)
\psline[showpoints=false,  linecolor =red]{-}(3,1)(4,2)
\uput*{.3}[150](4,2){a$_{3}$}

\psline[showpoints=true,  linecolor =red]{-}(4,2)(5,1)

\psline[showpoints=true,  linecolor =red]{-}(5,1)(6,0)

\psline[showpoints=false, linecolor = green]{-}(0,0)(1,-1)
\uput*{.3}[225](1,-1){b$_{1}$}

\psline[showpoints=false, linecolor = green]{-}(1,-1) (2, -2)
\uput*{.3}[225](2,-2){b$_{2}$}

\psline[showpoints=false, linecolor = green]{-}(2, -2)(3, -1)

\psline[showpoints=false, linecolor = green]{-}(3,-1) (4, 0)

\psline[showpoints=false, linecolor = green]{-}(4,0) (5, -1)
\uput*{.3}[225](5,-1){b$_{3}$}
\psline[showpoints=false, linecolor = green]{-}(5,-1) (6, 0)

\psline[showpoints =false, linecolor= cyan]{-}(2,2)(2,-3.6)
 \psline[showpoints =false, linecolor= cyan]{-}(4,0)(4,-3.6)
 \uput*{.2}[d](2,-3.6){$t_V(a_2)$}
  \uput*{.2}[d](4,-3.6){$t_W(a_2)$}
\psline[showpoints =false, linecolor= cyan]{<->}(2,-3.4)(4,-3.4)
\uput*{0}[r](2.9, -3.4){\textcolor{cyan}{t}}

 \psline[showpoints =false, linecolor= cyan]{-}(2,2)(8,2)
 \psline[showpoints =false, linecolor= cyan]{-}(4,0)(8,0)
\psline[showpoints =false, linecolor= cyan]{<->}(7.5, 2)(7.5,0)
\uput*{.2}[r](7.6, 0){$W(t_W(a_2))$}
\uput*{.2}[r](7.6, 2){$V(t_V(a_2))$}
\uput*{0}[d](7.7,1.1){\textcolor{cyan}{t}}

\psdots[linecolor=red](1,1)
\psdots[linecolor=red](2,2)
\psdots[linecolor=red](3,1)
\psdots[linecolor=red](4,2)
\psdots[linecolor=red](5,1)
\psdots[linecolor=red](6,0)

\psdots[linecolor=green](1,-1)
\psdots[linecolor=green](2,-2)
\psdots[linecolor=green](3,-1)
\psdots[linecolor=green](4, 0)
\psdots[linecolor=green](5,-1)

\end{pspicture}
\end{minipage}
\begin{center}
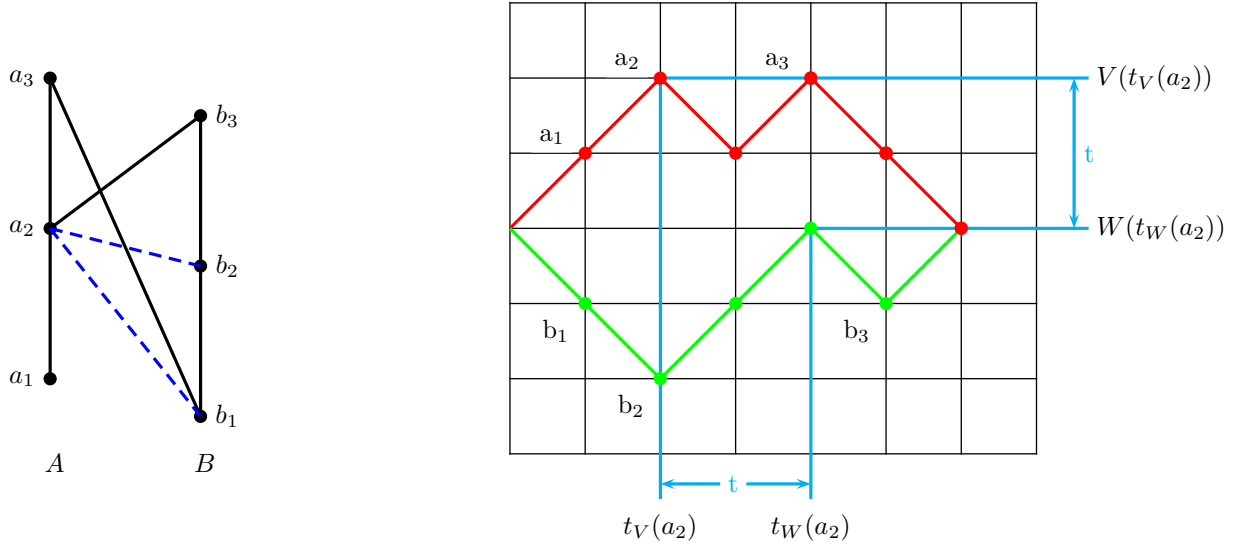
Figure 3: \footnotesize The description $I_C(a_2)$ in terms of the pair of walks.
\end{center}

\noindent time $t_W(i)$, the set
$\{\delta_1,\cdots,\delta_{t_W(i)}\}$ consists of $a_i$, the elements
below $a_i$ or incomparable to it, and no elements
that are above $a_i$.
Hence, $t_W(i) - t_V(i)$ is exactly the number of elements that
are incomparable to $a_i$. See for example, Figure 3.

Next we consider the second description of the incomparability window
for $a_i \in A$.
By part (a),
\begin{eqnarray*}
I_C(a_i) & = &  \ t_W(i) - t_V(i)\\
 & = & \ V(t_V(i)) - W(t_W(i))\\
& = & \ H_{V,W}(t_V(i)) + W(t_V(i)) - W(t_W(i))\\
& = & \ H_{V,W}(t_V(i)) + W(t_V(i)) - W(t_V(i) + I_T(a_i))\\
\end{eqnarray*}
The second equality above follows because at times $t_V(i)$ and
$t_W(i)$ respectively, both $V$ and $W$ have taken $i$ steps up and
$V(0)=W(0)$ (see for example Figure 3). The last equality follows
from the first equality when we substitute for $t_W(i)$.
Equation (\ref{eq:size2}) follows by a symmetric argument.
\end{proof}

Via linear interpolation, we may view $(V, W)$ as a pair of polygonal
paths in $\mathbb R^2$.
Let $A(C)$ denote the area between the two (piecewise linear) curves $V, W$.
As a consequence of the preceding identification, we note that the sum of the
incomparability windows can be written in terms of the area between
the curves:
\begin{corollary}\label{cor:window-is-area} $\displaystyle\sum_{c
    \in [n]} I_C(c) =  A(C).$
\end{corollary}
\begin{proof}
We will show that $\sum_{i = 1 }^k I_C(a_i) = \frac{1}{2}A(C) =
  \sum_{i = 1}^{n-k} I_C(b_i).$
The first equality follows by noticing that the lengths
  of the intervals $t_W(i) - t_V(i)$ are exactly the number of squares
  on the $i$-th diagonal (running from top left to bottom right)
  between $V$ and $W$ (see Figure 1). Since
  each square has a
  side length of
  $\sqrt{2}$, the area is 2.
The second equality follows by counting the squares
between $V$ and $W$ from bottom left to top right.
\end{proof}
\begin{proof}[Proof of Theorem \ref{thm:mapping-phi}]
Fix a poset $P
  \in \mathcal P_1$ and let $(V,W) \in \Phi(P)$. Let $C =
  \Gamma^{-1}(V,W) = (A,B,\prec)$ be the two-chain cover corresponding to
  $(V,W)$. Recall from
  (\ref{eq:tau}) that for an element $x$, $\tau_P(x)$ is the number of
  elements of the poset smaller than or equal to $x$.
Recall that we can think of the elements of $A$ being labeled $a_1,
  \cdots,a_k$ and those of $B$ being labeled $b_1, \cdots,b_{n-k}$.
Let $\bar C$ be the unlabeled two-chain cover of which $C$ is the canonical
representative. Since
$C \in \Psi(P)$,  for any labeled two-chain cover $C_P$ of $P$
in $\bar C$, there is a two-chain isomorphism $\ga$ between  $C_P$ and
$C$. Consequently, for any element $x \in [n]$,
\begin{align*}
I_P(x) = I_{C_P}(x) = I_C(\ga(x)).
\end{align*}
Also, if $\alpha(x) = a_i \in A$ then $\tau_P(x) = t_V(i) $ while if
$\alpha(x) = b_j \in B$, then $\tau_P(x) = s_W(j)$.
Therefore, using the bounds (\ref{eq:size}) and (\ref{eq:size2}) from Lemma
\ref{lem:window-is-intercept}, we have
\begin{align*}
|I_{P}(x)- & H_{V,W}(\tau_P(x))| = |I_C(\alpha(x) ) -
 H_{V,W}(\tau_P(x))| \\
 \leq &  \ |V(\tau_P(x) + I_C(\alpha(x))) - V(\tau_P(x))| + |W(\tau_P(x) +
I_C(\alpha(x))) -
 W(\tau_P(x))| \\
= & \ |V(\tau_P(x) + I_P(x)) - V(\tau_P(x))| + |W(\tau_P(x) + I_P(x)) -
 W(\tau_P(x))|.
\end{align*}
\end{proof}

Let $(V,W)$ be a random pair of non-hitting walks drawn from
$\BbbP_\texttt{cov}$
and let $C = (A,B,\prec~)$ be the
corresponding unlabeled two-chain cover. Let $U_n$ be uniform over
$[n]$. We show that
$H_{V,W}(\tau_{V,W}(U_n))$ is distributionally
invariant if we replace $\tau_{V,W}(U_n)$ with an independent
uniform variable $U_n'$ on $\{1, \cdots,n\}$.
\begin{proof}[Proof of Theorem \ref{Lsym}]

Let  $m \in 2\mathbb{Z}_{\geq}$ be fixed and define
\begin{equation}
\mathcal A_{t}(e_1, e_2)  :=  \{ (V,W): H_{V,W}(t) = m,\: (\pi_t(V),
\pi_t(W)) = ( e_1, e_2)\}
\end{equation}
where $ e_i \in \{-1, +1\}.$
Define an involution $\imath : \BB_n \to  \BB_n $ as follows:
\begin{align*}
\imath(V,W) := (V', W') \text{ where:}
\begin{cases}
V'(0)  = W'(0) =0,  \\
\pi_i(V') = - \pi_{n-i+1}(V),\\
  \pi_i(W') = - \pi_{n-i+1}(W)
  \end{cases}
\end{align*}
 \text{ for all $ 1 \le i \le n$}.  In words $\imath$ maps $(V, W)$ to
 its time reversal, additionally shifting the height so that $V(n),
 W(n)$ maps to $(0, 0)$.
 It is easy to see that $ (V', W') \in \BB_n$ and that $ \imath^2
 = id$.  This implies that $\imath$ is a bijection on $\BB_n$
and thus preserves the uniform distribution $\BbbP_\texttt{cov}$ over
 $\BB_n$.

We claim that
\begin{align}\label{eq:reversal}
\sum_{t=1}^n |\mathcal A_{t}(-1, +1)|  =  \sum_{t=1}^n |\mathcal
 A_{t}(+1, -1)| .
\end{align}
This follows by observing that
$\imath$ maps $\mathcal A_{t}(-1, +1)$ into  $\mathcal A_{n-t+1}(+1,
-1)$ and is
a  bijection.
Indeed: $
H_{\imath(V,W)}(n-t+1) = H_{V,W}(t)$ and if $(\pi_k(V),\pi_k(W))  = (-1,
+1)$, then $(\pi_{n-t+1}(V'),\pi_{n-t+1}(W'))  = (+1, -1)$ as well.

Let $C = (A,B,\prec)$ be the two-chain cover corresponding to $(V,W)$.
Suppose for $c \in [n]$, $\tau_{V,W}(c) = t$. Then $\pi_{t}(V)
= +1$ if and only if $c \in A$. Moreover there is exactly
one element of $A$ such that $\tau_{V,W}(c) = t$. Similarly,
$\pi_{t}(W)
= -1$ if and only if $c \in B$, and there is exactly
one such element $c$ in $B$.
Therefore, for $1\leq t \leq n$, we have
\begin{align}\label{eq:nt}
|\{c \in [n]:
\tau_{V,W}(c) = t\}| = \left\{ \begin{array}{lll}
1 &  \text{ if} & \pi_t(V) = +1, \pi_t(W) = +1 \ \\
1 &  \text{ if} & \pi_t(V) = -1, \pi_t(W) = -1 \ \\
2 &  \text{ if} & \pi_t(V) = +1, \pi_t(W) = -1 \ \\
0 &  \text{ if} & \pi_t(V) = -1, \pi_t(W) = +1. \ \\
\end{array} \right.
\end{align}

Therefore, for every $m \in 2 \mathbb Z_{\geq}$,
\begin{align*}
\mathbb P_{(V,W) \sim \BbbP_\texttt{cov}}(H_{V,W} & (\tau_{V,W} (U_n))=m)
\\
= & \
|\BB_n|^{-1} n^{-1} \sum_{(V,W) \in \BB_n} \sum_{c \in [n]} \mathbb
I(H_{V,W}(\tau_{V,W}(c)) = m) \\
= & \
|\BB_n|^{-1} n^{-1} \sum_{(V,W) \in \BB_n}
\sum_{t=1}^n |\{c \in [n]:
\tau_{V,W}(c) = t\}| \mathbb I(H_{V,W}(t) = m). \\
\end{align*}
Switching the order of summation and using (\ref{eq:nt}), we have
\begin{align*}
\mathbb P_{(V,W) \sim \BbbP_\texttt{cov}}(H_{V,W} & (\tau_{V,W} (U_n))=m)\\
= & \ |\BB_n|^{-1}
n^{-1} \sum_{t=1}^n |\mathcal
A_{t}(+1, +1)| + |\mathcal A_{t}(-1, -1)| + 2 |\mathcal A_{t}(+1,
-1)| \\
= & \ |\BB_n|^{-1} n^{-1} \sum_{t=1}^n \sum_{\{e_1, e_2\}}
|\mathcal A_{t}(e_1, e_2)| \mathrm{\  \ \ \ \ \ \  \ \ \ \ \ \ \ \ \
  (by \ (\ref{eq:reversal}))} \\
= & \ |\BB_n|^{-1} n^{-1} \times
\sum_{t = 1}^n | \{ (V,W): H_{V,W}(t) = m\} |\\
= & \ \mathbb P(H_{V,W}(U_n') =m).
\end{align*}
\end{proof}
\section{Proof of Theorem
  \ref{thm:SMain}}\label{sec:convergence-for-one-factor}~In this
section we prove that the size of the incomparability window
of a typical
element for a random poset with one factor converges to a Brownian
excursion at a uniformly random time under the correct scaling.
\begin{theorem}
  \label{thm:SMain}
If $P$ is  chosen uniformly at random from $\PP_1$ and $U_n$ is
uniform over the ground set $[n]$ independent of $P$, then
 \begin{equation}
 \label{Eq:Main}
 \frac{I_{P}(U_n)}{\sqrt n}  \stackrel{D}{ \Rightarrow} \frac{1}{\sqrt
 2} B^{\text{\emph{ex}}}(U),
 \end{equation}
  where $B^{\text{\emph{ex}}}$ is standard Brownian Excursion on
  $[0,1]$ and $U \sim$ Uniform $[0,1]$ independent of
  $B^{\text{\emph{ex}}}$.
  \end{theorem}
Recall that by Lemma \ref{lem:reduction-to-one-factor}, this implies the main
result.
\subsection{The Distribution of $\BbbP_\texttt{cov}$}
Instead of working
directly with random posets we will work with the corresponding random
walks as described in Section
\ref{sec:combinatorial-maps} and hence
we begin with a convenient description of the uniform measure on $\BB_n$. Recall that by Lemma \ref{lem:sig}, $\BbbP_\texttt{cov}$ is the uniform measure on $\BB_n$.
Let $Z^{(1)}$ and $Z^{(2)}$ be two independent random walks starting
at the origin with
\begin{align}Z_k^{(i)} = \sum_{l=1}^k z_l^{(i)},
\end{align}
where the $z_l^{(i)}$
are \iid \  $\pm 1$ with equal probability for $i=1,2$. Let $(V,
W) $ be chosen uniformly at random from $\BB_n$. Let
$X \stackrel{D}{=} Y \ | \
Z$ denote that the random variable $X$ is
distributed like $Y$ conditioned on the event $Z$.
By definition,
\begin{align}\label{eq:dist}
(V, W) \stackrel{D}{=} ( Z^{(1)}, Z^{(2)} )  \ | \ Z_k^{(1)}
>  Z_k^{(2)} \text{ for } 0< k < n, Z_n^{(1)}= Z_n^{(2)} .
\end{align}

Consider the random walks $ Z^{(1)} +Z^{(2)}$  and $Z^{(1)}
-Z^{(2)} $. Note that unlike $ Z^{(1)}$ and $Z^{(2)}$, these random
walks are dependent. Their joint distribution can be described in the
following way:

Let $ \{ \xi_i\}, \{ \xi'_i\}$ for $1 \leq i \leq n$ be \iid \ random
variables taking values $ \pm 1$ with probability $ 1/2$. Let $ u_i$
be \iid\ $\{0,1\}$ independent of $\{ \xi_i\} , \{ \xi'_i\}$.
If we think of the variable $u_i$ being 1 if $z_i^{(1)} =
z_i^{(2)}$ and 0 otherwise then it can be seen that $\frac{1}{2}
(z_i^{(1)}+ z_i^{(2)})$ and $\frac{1}{2}
(z_i^{(1)}- z_i^{(2)})$
are distributed
according to $\xi_i u_i$ and $\xi'_i (1-u_i)$ respectively.
For $1
\leq k \leq n$, let
\begin{align}\label{eq:y-y'}
Y_k & := \xi_1 u_1 +\xi_2 u_2 + \ldots+ \xi_k u_k \\
Y'_k & := \xi'_1(1- u_1) +\xi'_2(1- u_2) + \ldots+ \xi'_k (1-u_k).
\end{align}
Then,
\begin{align}\label{eq:z-to-y}
\left (   \frac{1}{2}(Z^{(1)} +Z^{(2)} ),  \frac{1}{2}(Z^{(1)} -Z^{(2)} )
\right ) \stackrel{D}{=}  ( Y , Y').
\end{align}
Thus, we have the following
distributional identity.
\begin{lemma}\label{lem:dist_eq}
Let $(V,W)$ be uniform on $\BB_n$ and let $Y,Y'$ be as
above. Then
\begin{equation}\label{eq:jointwalk}
\left( \frac{1}{2}(V + W), \frac{1}{2}(V - W) \right)
\stackrel{D}{=} ( Y , Y') \ \big | \ Y'_k > 0 \text{ \ for \ }0 < k < n, Y'_n
=0.
\end{equation}
\end{lemma}
\begin{proof}
\begin{align*}
\left( \frac{1}{2}(V + W), \frac{1}{2}(V - W) \right)
& \stackrel{D}{=}
\left( \frac{1}{2}(Z^{(1)}+Z^{(2)}),
  \frac{1}{2}(Z^{(1)}-Z^{(2)})  \right)
\ | \ Z \ \ \ \ \ \ \ \ \ \ \ \text{(by \eqref{eq:dist})}\\
& \stackrel{D}{=}
( Y , Y') \ \big | \ Y'_k > 0 \text{ for }0 < k < n, Y'_n
=0  \ \ \ \ \ \ \ \ \  \ \ \ \text{(by \eqref{eq:z-to-y})}\\
\end{align*}
where above $Z$ is the event
\begin{align*}
\{ Z_k^{(1)}- Z_k^{(2)} > 0 \text{ for } 0< k < n, Z_n^{(1)}-
Z_n^{(2)} = 0\}.
\end{align*}
\end{proof}
As an immediate corollary we obtain the distribution of the second marginal:
\begin{corollary}
\label{cor:Excursion}
Let $Y'_k = x_1 + x_2 + \ldots + x_k$ where $\{ x_i \}_{ 1 \le i \le n}$
are \iid \ with $ \bP(x_1 =0) =1/2$, $\bP(x_1 =+1) = \bP(x_1 =-1) =
1/4$. Let $(V, W)$ be chosen uniformly from $\BB_n$. Then
\begin{equation}
\frac{1}{2}(V - W) \stackrel{D}{=}Y' \ \big | \ Y'_k > 0  \text{ for
}0 < k < n, Y'_n =0.
\end{equation}
\end{corollary}

\subsection{An Upper Bound on the incomparability Window
    $I_{P}(U_n)$}
Corollary \ref{cor:Excursion} allows us to find the limiting
distribution of the
``average'' incomparability window.  For a poset $P \in \PP_1$, define
\begin{equation}
\hat I \left(P\right) = \frac{1}{n} \sum_{x \in [n]} I_{P}(x)
\end{equation}
\begin{theorem}\label{thm:convavg}
If $P$ is  chosen uniformly at random from $\mP_0$, then
\[ \frac{\hat I(P)}{\sqrt{n}} \stackrel{D}{ \Rightarrow} \sqrt{2} \int_0^1
B^{\text{\emph{ex}}}(t) dt, \]
where $B^{\text{\emph{ex}}}(\cdot)$ is a standard Brownian excursion
on $[0,1]$.
  \end{theorem}
 \begin{proof}
 By Corollary \ref{cor:window-is-area}, $\hat I(P) =
 n^{-1} A(C)$ where for any two-chain cover $C$ of $P$,
 $A(C)$ denotes the area between
 the piecewise-linear interpolations of the corresponding walks $V$
 and $W$ over the time
 interval $[0, n]$.
By abuse of notation, denote the pair of piecewise linear
 interpolated curves by $V_s, W_s$ for $s \in [0, n]$.  Then
\begin{equation}
\frac{1}{n} A(C) = \frac{1}{n} \int_{0}^n (V_s -
W_s) \: ds = \sqrt{n} \int_{0}^1 \frac{1}{\sqrt n}
\left(V_{n t} - W_{n t }\right) dt.
\end{equation}

Now we apply Corollary \ref{cor:Excursion}.  Let $\gs^2 := \E x_i^2 = 1/2.$
From the  invariance principle for random walks excursions ( see
\cite{K76} for example),
\begin{equation}\label{eq:convergence-to-bex}
\left (\frac{1}{\gs \sqrt{ n }} Y'_{nt}, 0 \le t\le 1
  \Big |Y'_k > 0  \text{ for }0 < k < n, Y'_n =0 \right )
\stackrel{D}{ \Rightarrow} \left ( B^{\text{ex}}(t), 0 \le t\le 1
\right ),
\end{equation}
where $Y'_{nt}$ is the piecewise linear interpolation of $Y'$.
Here the weak convergence takes place in $C[0,1]$ equipped
with uniform topology.
It follows from the continuous mapping theorem that
\begin{equation}
\int_{0}^1 2 \frac{1}{\sqrt{n}} Y'_{n t } dt
\stackrel{D}{ \Rightarrow} \int_0^1 \sqrt{2} B^{\text{\emph{ex}}}(t) dt.
\end{equation}
The lemma now follows.
\end{proof}
\begin{corollary}\label{cor:max}
Let $P \in \PP_1$ be uniformly distributed and let $U_n$ be chosen
uniformly from the ground set $[n]$ independent from $P$. Then
\begin{equation}
\pr (I_{P}(U_n) > n^{2/3}) \rightarrow 0
\end{equation}
\end{corollary}

\begin{proof}
From Theorem \ref{thm:convavg} it follows that $\mathbb P
(\hat{I}(P) \leq n^{1/2+1/12}) \to 1 $ as $n \to
\infty$. Markov's inequality implies that
$\mathbb P(I_P(U_n) > n^{2/3} | \hat{I}(P) \leq n^{1/2+1/12}) \to 0$
as $n \to \infty$ and the proof follows.
\end{proof}

\subsection{Conclusion of the Proof of Theorem \ref{thm:SMain}}
In this section we show that for a randomly
chosen poset $P \in
\PP_1$, and for $U_n $ a uniformly random element of $[n]$,
$$Err(P, U_n):= |V(\tau_P(U_n) + I_P(U_n)) -
V(\tau_P(U_n))| +
 |W(\tau_P(U_n) + I_P(U_n)) -
 W(\tau_P(U_n))|$$ behaves as $o(\sqrt n)$ in probability. We
show that for typical elements of the poset, the fluctuations in the
corresponding random walk representation are small.

Let $Y_k$ and $Y_k'$ be the random walks as defined in \eqref{eq:y-y'}.
The next result shows that the probability that either $Y_k$ or $Y_k'$ exhibits a deviation of more than
$n^{1/3+\delta}$ in
$n^{2/3}$ steps of the walks goes to 0. Further, the probabilities go
to zero even when conditioned on the event $\{Y'_k > 0 \ \forall \: 1 \le
k \le n-1, \ Y'_n =0 \} $.
Together with \eqref{eq:jointwalk} the result will imply that the
walks $V,W$ have  fluctuations bounded by $n^{1/3 + \delta}$ if we look
at most $n^{2/3}$ steps away, with high probability.
We make use of the following result for a simple symmetric random walk (SSRW).
\begin{proposition}[\cite{Fe68},
  Lemma 2, page 78] \label{prop:first-return-to-zero}
For every $k \in \mathbb N$,
\begin{align}
\pr( \text{ The first return to $0$ of a SSRW is at } 2k  ) =
\frac{1}{2k-1}{2k  \choose k} 2^{-2k}  = \Theta \left(\frac{1}{k^{3/2}}\right).
\end{align}
\end{proposition}

\begin{lemma}\label{fluctuations}
Let $Y, Y'$ be as in \eqref{eq:y-y'} and let $G = \{Y'_k > 0 \
\forall 1 \le k \le n-1, \
Y'_n =0 \} $. Then, for all $1 \leq k \leq n$ and for any $\gd> 0$,
\begin{align}\label{eq:fluc1}
\pr \left ( \max_{ 0 \le s \le n^{2/3} }
  |Y_{k+s} - Y_k| > n^{1/3 + \gd} \Big | G\right) \rightarrow 0
\end{align}
and
\begin{align}\label{eq:fluc2}
\pr \left ( \max_{ 0 \le s \le n^{2/3} }
  |Y'_{k+s} - Y'_k| > n^{1/3 + \gd} \Big | G\right) \rightarrow 0
\end{align}
as $n \rightarrow \infty$.
\end{lemma}
\begin{proof}
For any $1 \leq k \leq n$,
we can bound the probability of the unconditional events above by the
Chernoff bound as follows
\begin{align}\label{eq:chernoff-for-deviation-y}
\nonumber \pr \left( \max_{ 0 \le s \le n^{2/3} }
  |Y_{k+s} - Y_k| > n^{1/3 + \gd} \right)
\leq & \sum_{0 \leq s \leq n^{2/3}} \pr \left(
  |Y_s| > n^{1/3 + \gd} \right) \\
\leq  & \ 2 n^{2/3} e^{-n^{2\gd/3}/4}.
\end{align}
Similarly,
\begin{align}\label{eq:chernoff-for-deviation-y'}
\pr \left( \max_{ 0 \le s \le n^{2/3} }
  |Y'_{k+s} - Y'_k| > n^{1/3 + \gd} \right)
\leq  & \ 2 n^{2/3} e^{-n^{2\gd/3}/4} .
\end{align}
Further, we have
\begin{align*}
& \ \pr( \text{ The first return to $0$ of $Y'_k$ is at } n  ) \\
= & \  \sum_{0 \leq i \leq n}^{*} \pr( |\{j: u_j = 1)\}| = i) \pr (\text{
  the first return to $0$ of a SSRW is at } n-i ) \\
 = & \  \Theta \left( \frac{1}{n^{3/2}} \right) \  \sum_{0 \leq i \leq
  n}^{*}
  \pr( |\{j: u_j = 1)\}|
= i)  \ \ \ \ \ \ \ \ \ \ \text{(by Proposition
  \ref{prop:first-return-to-zero} )} \\
\end{align*}
where the sum $\sum^*$ is restricted to those $i$ for which $n-i \in 2
\mathbb Z_{\geq}$,
and therefore
\begin{align}\label{eq:positive-excursions}
\pr & \left( G \right) =
\Theta \left( \frac{1}{n^{3/2}} \right).
\end{align}
By \eqref{eq:chernoff-for-deviation-y} and
\eqref{eq:positive-excursions}, we have
\begin{align*}
\pr \left( \max_{ 0 \le s \le n^{2/3} }
  |Y_{k+s} - Y_k| >  n^{1/3 + \gd} \Big | G \right)
\ \leq & \ \frac{\pr \left( \displaystyle\max_{
      0 \le s \le
      n^{2/3} }
  |Y_{k+s} - Y_k| >  n^{1/3 + \gd} \right)}{\pr (G) } \\
& \leq O(1)\frac{n^{2/3}e^{-n^{2 \gd / 3}/4}}{n^{-3/2}}
\ \rightarrow 0.
\end{align*}
A similar argument using \eqref{eq:chernoff-for-deviation-y'} and
\eqref{eq:positive-excursions} shows that
\begin{align*}
& \ \pr \left( \max_{ 0 \le s \le n^{2/3} }
  |Y'_{k+s} - Y'_k| >  n^{1/3 + \gd} \Big | G \right) \rightarrow 0.
\end{align*}
\end{proof}

\begin{corollary} \label{cor:fluc}
If $( V, W)$ are drawn uniformly from $\BB_n$, then for all $1 \leq k
\leq n$ and $\delta > 0$,
\begin{equation*}
\pr \left ( \max_{ 0 \le s \le  n^{2/3} }
  |V(k+s) - V(k)| > 2n^{1/3 + \gd}  \right) \rightarrow 0
\end{equation*}
and
\begin{equation*}
\pr \left ( \max_{ 0 \leq s \le n^{2/3} }
  |W(k+s) - W(k)| > 2n^{1/3 + \gd}  \right) \rightarrow 0 .
\end{equation*}
\end{corollary}

\begin{proof}
The claim follows by adding the bounds from \eqref{eq:fluc1} and
\eqref{eq:fluc2} from Lemma \ref{fluctuations} and the description of
$(V, W)$ in (\ref{eq:jointwalk}).
\end{proof}

We can now complete the proof of Theorem \ref{thm:SMain}.
\begin{proof}[Proof of Theorem \ref{thm:SMain}]
Recall that we would like to show that
\begin{equation}
 \frac{1}{\sqrt{n}} I_{P}(U_n) \stackrel{D}{ \Rightarrow}
 \frac{1}{\sqrt{2}}B^{\text{\emph{ex}}}(U).
 \end{equation}
Let $(P, (V, W))$  be chosen according to $2$-dimensional marginal of
$\BbbP$ on $\PP_1 \times \BB_n$ and  let $U_n$
be chosen uniformly at random from $[n]$ independent of $(P,
(V, W))$. By Corollary \ref{cor:max}, with
probability going to 1, $I_{P}(U_n) \leq n^{2/3}$. Therefore,
by part (2) of Theorem \ref{thm:mapping-phi}, with high probability
\begin{align}
\nonumber |I_{P}(U_n) - H_{V,W}(\tau_{P}(U_n))| \le
\max_{0 \le k \le n^{2/3}} |V(\tau_P(U_n) + k) -
V(\tau_P(U_n))| \\
+
\max_{0 \le
  k \le n^{2/3}} |W( \tau_P ( U_n)+ k) - W(\tau_P(U_n))|.
\end{align}
By Corollary
\ref{cor:fluc},
\[\frac{1}{\sqrt n}\left( \max_{0 \le k \le n^{2/3}}
  |V(\tau(U_n) + k) -
  V(\tau(U_n))| + \max_{0 \le
  k \le n^{2/3}} |W( \tau
  (U_n)+
  k) - W(\tau(U_n))|\right) \stackrel{P}{\to} 0
, \]
where $ \stackrel{P}{\to}$ denotes convergence in probability.
Therefore,
\begin{align}
 \frac{1}{\sqrt n}|I_{P}(U_n) - H_{V,W}(\tau_{P}(U_n))|
\stackrel{P}{\rightarrow} 0.
\end{align}
By Lemma \ref{lem:tauP-tauT}, $H_{V,W}(\tau_P(U_n))$ equals
$H_{V,W}(\tau_{V,W}(U_n'))$ where $U_n'$ is some
uniform random variable on $[n]$ which is independent of $(V, W)$.
Hence,
\begin{align}
\nonumber \frac{1}{\sqrt n} |I_{P}(U_n) -
H_{V,W}(\tau_{V,W}(U_n'))|
\stackrel{P}{\rightarrow} 0.
\end{align}

Let $(V,W)$ be chosen according to $\BbbP_\texttt{cov}$.
The invariance principle for
random walk excursions implies (see \ref{eq:convergence-to-bex}) that
\[ \frac{1}{\sqrt n} H_{V,W}(U'_n) \stackrel{D}{ \Rightarrow}
\frac{1}{\sqrt 2} B^{\text{ex}}(U).\]
By Theorem \ref{Lsym},  $H_{V,W}(\tau_{V,W} (U_n'))$ has
same distribution as $H_{V,W}(U'_n)$ and therefore
\[ \frac{1}{\sqrt n} H_{V,W}(\tau_{V,W}(U_n'))
\stackrel{D}{ \Rightarrow}
\frac{1}{\sqrt 2} B^{\text{ex}}(U). \]
This implies the claim
\[ \frac{1}{\sqrt n} I_P(U_n)
\stackrel{D}{ \Rightarrow}
\frac{1}{\sqrt 2} B^{\text{ex}}(U) .\]
\end{proof}

\section{Distribution of Height}
In this section, we prove Theorem \ref{thm:height} showing that the
height of random width-2 poset is $n/2$ with a Gaussian fluctuation.
\begin{proof}[Proof of Theorem \ref{thm:height}]
Let $M_1\prec M_2\prec \ldots \prec M_k$ be the factors
of $P$ with 
$C_i = (A_i, B_i, \prec)$ being a two-chain cover of $M_i$. The
longest chain in $P$ can be constructed by concatenating the
longer of the two chains from each of the $k$ factors. Thus the height
of $P$ is given by
\[ h(P) = \sum_{i=1}^k \max( |A_i|, |B_i|).\]
As in Lemma \ref{lem:BG-large-factor}, it is easy to check that
that the height is essentially determined by the largest factor which
is of size at least $n-O(\log n)$. We may
proceed analogously to the proof of Theorem
\ref{thm:distribution-of-window} to conclude that it is enough to
prove the above theorem for the special case when the
poset $P$ is chosen uniformly at random  from $\PP_1$, i.e. it has
only one factor.

For $P \in \PP_1$, let $(V, W) \in \Phi(P)$ be an associated
pair of walks. It is easy to see that
\[ h(P) = \frac{n+ |V(n)|}{2} = \frac{n+ |W(n)|}{2}. \]

Thus, we would like to find the limiting distribution of $(V(n) +
W(n)) / 4 \sqrt n$. The next lemma completes the proof once we
recall the distributional identity given in Lemma \ref{lem:dist_eq}.
\end{proof}

\begin{lemma}
Let $Y, Y'$  be as in  the paragraph preceding Lemma \ref{lem:dist_eq}. Then
\[ \frac{2Y_n}{\sqrt n } \big | \ Y'_k > 0 \text{ \ for }0 < k < n, Y'_n
=0 \stackrel{D}{\Rightarrow} N(0,1)\]
as $n \to \infty$.
\end{lemma}
\begin{proof}
Let $\Phi(\cdot)$ denote the standard normal distribution
function. Let $G$ denote the event  $ \{Y'_k > 0 \text{ for }0 < k <
n, Y'_n =0 \}$. Fix $x \in \mathbb R$. We then have
\begin{align*}
\pr\left( \frac{2Y_n}{\sqrt n }  \le x | G \right)  - \Phi(x) &= \sum^*_{0 \le m
  \le n} \left( \pr\left( \frac{2Y_n}{\sqrt n }  \le x  | G, \sum_{i=1}^n
  u_i = m\right)  - \Phi(x)\right) \pr\left(  \sum_{i=1}^n u_i = m |
  G\right) \\
&= \sum^*_{0 \le m \le n} \left( \pr\left( \frac{2}{\sqrt n } \sum_{i=1}^m
  \xi_i \le x \right)  - \Phi(x)\right) \pr\left(  \sum_{i=1}^n u_i = m | G\right)
\end{align*}
where the sum $\sum^*$ always includes the added restriction $n-m \in
2 \mathbb Z$. The lemma now follows from a simple application of the
Berry-Ess\'{e}en bound (see for example \cite{Dur}) once we prove that
$\sum_{i=1}^n u_i$ is
concentrated around $n/2$ even when conditioned on $G$. More
precisely, we want to show,
\[ \pr\left( |\sum_{i=1}^n u_i - n/2| > n^{3/4} |G \right) \to 0.\]
Note that \[ \pr\left( |\sum_{i=1}^n u_i - n/2| > n^{3/4} |G \right) \le
\displaystyle \frac{\pr\left( |\sum_{i=1}^n u_i - n/2| > n^{3/4}  \right)
}{\pr(G)} \le \frac{c_1 \exp( - c_2 n^{1/2} )}{\pr(G)},\]
for suitable constants $c_1, c_2 > 0$.
By \eqref{eq:positive-excursions}
\begin{align*}
\pr(G)   = \Theta \left(\frac{1}{n^{3/2}}\right) ,
\end{align*}
and the claim follows.
\end{proof}

\section{Conclusions}\label{sec:conclusions}
The results in this paper should be viewed as a first step in the
analysis of random posets of bounded width. Our results provide the
asymptotic distribution of the number of elements incomparable to a
random element. However, more detailed information is desirable.

In our results we find a distribution $F$ such that for a random element
$x$, there is a maximal chain of
 incomparable elements $y_{i(0)} \prec \ldots \prec y_{j(0)}$ and $j(0)-i(0)$ is
distributed according to $F$.
 In fact one would expect to extract more detailed information about
the ``neighborhood" of $x$ and that it has the following structure:
 $x = x_0$ belongs to a chain $\prec \cdots \prec x_{-1} \prec x_0 \prec x_1 \prec \cdots$
 and there
exists another chain $\cdots \prec y_{-1} \prec y_0 \prec y_1 \prec \cdots$
 such that each element $x_k$ of the first chain is incomparable to elements
 $y_{i(k)} \prec \cdots \prec y_{j(k)}$ of the other chain. Then it would be desirable
to identify the (joint)
 distribution of $i(k+1)-i(k)$ and $j(k+1)-j(k)$ for small values
 of $k$. Even more detailed information is desirable in terms of the joint
distribution of the $x$'s, the $y$'s and the elements incomparable to
$y$'s etc.
As mentioned in the introduction, it would also be desirable to extend
the analysis here to posets of bounded width greater than $2$.

\section{Acknowledgments}

E.M. is grateful to Graham Brightwell for discussions which initiated
this work. The authors are grateful to David Aldous and Steve Evans
for  pointing out the distribution of the height of a Brownian
excursion.

\end{document}